\documentclass[10pt]{article}

\usepackage{amsmath,amssymb,amsthm}
\usepackage[left=3cm,right=3cm,top=2.9cm,bottom=2.9cm]{geometry}

\usepackage{color,enumitem,graphicx}
\usepackage[colorlinks=true,urlcolor=blue,
citecolor=red,linkcolor=blue,linktocpage,pdfpagelabels,
bookmarksnumbered,bookmarksopen]{hyperref}
\usepackage[hyperpageref]{backref}

\newcommand{\bdry}[1]{\partial #1}
\newcommand{\bgset}[1]{\big\{#1\big\}}
\newcommand{\C}{{\cal C}}
\newcommand{\F}{{\cal F}}
\newcommand{\calN}{{\cal N}}
\newcommand{\eps}{\varepsilon}
\newcommand{\Fucik}{Fu\v c\'\i k }
\newcommand{\Il}[1]{\text{I}_{#1}}
\newcommand{\Iu}[1]{\text{I}^{#1}}
\newcommand{\II}[1]{\text{II}_{#1}}
\newcommand{\ip}[3][]{\left(#2,#3\right)_{#1}}
\newcommand{\isom}{\approx}
\newcommand{\M}{{\cal M}}
\newcommand{\norm}[2][]{\left\|#2\right\|_{#1}}
\renewcommand{\o}{\text{o}}
\newcommand{\pnorm}[2][]{\if #1'' \left|#2\right|_p \else \left|#2\right|_{#1} \fi}
\newcommand{\R}{\mathbb R}
\newcommand{\restr}[2]{\left.#1\right|_{#2}}
\newcommand{\seq}[1]{\left(#1\right)}
\newcommand{\set}[1]{\left\{#1\right\}}
\newcommand{\spnorm}[2][]{\if #1'' |#2|_p \else |#2|_{#1} \fi}
\newcommand{\wto}{\rightharpoonup}
\newcommand{\Z}{\mathbb Z}

\DeclareMathOperator{\divg}{div}

\newenvironment{enumalph}{\begin{enumerate}

}{\end{enumerate}}

\newenvironment{enumroman}{\begin{enumerate}

}{\end{enumerate}}

\newtheorem{proposition}{Proposition}[section]
\newtheorem{theorem}[proposition]{Theorem}

\numberwithin{equation}{section}

\title{\bf A note on the Dancer-\Fucik spectra of the fractional $p$-Laplacian and Laplacian operators\thanks{{\em MSC2010:} Primary 35R11, 35P30, Secondary 35A15
\newline \smallskip \indent\; {\em Key Words and Phrases:} fractional $p$-Laplacian, Dancer-\Fucik spectrum, critical groups}}
\author{\bf Kanishka Perera\\
Department of Mathematical Sciences\\
Florida Institute of Technology\\
Melbourne, FL 32901, USA\\
[\bigskipamount]
\bf Marco Squassina\thanks{The second-named author was supported by 2009 MIUR project: ``Variational and Topological Methods in the Study of Nonlinear Phenomena''.}\\
Dipartimento di Informatica\\
Universit\`a degli Studi di Verona\\
37134 Verona, Italy\\
[\bigskipamount]
\bf Yang Yang\thanks{This work was completed while the third-named author was visiting the Department of Mathematical Sciences at the Florida Institute of Technology, and she is grateful for the kind hospitality of the department. Project supported by NSFC-Tian Yuan Special Foundation (No. 11226116), Natural Science Foundation of Jiangsu Province of China for Young Scholars (No. BK2012109), and the China Scholarship Council (No. 201208320435).}\\
School of Science\\
Jiangnan University\\
Wuxi, 214122, China}
\date{}

\begin{document}

\maketitle

\begin{abstract}
We study the Dancer-\Fucik spectrum of the fractional $p$-Laplacian operator. We construct an unbounded sequence of decreasing curves in the spectrum using a suitable minimax scheme. For $p = 2$, we present a very accurate local analysis. We construct the minimal and maximal curves of the spectrum locally near the points where it intersects the main diagonal of the plane. We give a sufficient condition for the region between them to be nonempty, and show that it is free of the spectrum in the case of a simple eigenvalue. Finally we compute the critical groups in various regions separated by these curves. We compute them precisely in certain regions, and prove a shifting theorem that gives a finite-dimensional reduction in certain other regions. This allows us to obtain nontrivial solutions of perturbed problems with nonlinearities crossing a curve of the spectrum via a comparison of the critical groups at zero and infinity.
\end{abstract}

\section{Introduction}

For $p \in (1,\infty)$, $s \in (0,1)$, and $N > sp$, the fractional $p$-Laplacian $(- \Delta)_p^s$ is the nonlinear nonlocal operator defined on smooth functions by
\[
(- \Delta)_p^s\, u(x) = 2 \lim_{\eps \searrow 0} \int_{\R^N \setminus B_\eps(x)} \frac{|u(x) - u(y)|^{p-2}\, (u(x) - u(y))}{|x - y|^{N+sp}}\, dy, \quad x \in \R^N.
\]
This definition is consistent, up to a normalization constant depending on $N$ and $s$, with the usual definition of the linear fractional Laplacian operator $(- \Delta)^s$ when $p = 2$. There is currently a rapidly growing literature on problems involving these nonlocal operators. In particular, fractional $p$-eigenvalue problems have been studied in Lindgren and Lindqvist \cite{MR3148135}, Iannizzotto and Squassina \cite{MR3245079}, and Franzina and Palatucci \cite{FrPa}, regularity of fractional $p$-minimizers in Di Castro et al.\! \cite{DiKuPa}, and existence via Morse theory in Iannizzotto et al.\! \cite{IaLiPeSq}. We refer to Caffarelli \cite{Ca} for the motivations that have lead to their study.

Let $\Omega$ be a bounded domain in $\R^N$ with Lipschitz boundary $\bdry{\Omega}$. The Dancer-\Fucik spectrum of the operator $(- \Delta)_p^s$ in $\Omega$ is the set $\Sigma_p^s(\Omega)$ of all points $(a,b) \in \R^2$ such that the problem
\begin{equation} \label{1}
\left\{\begin{aligned}
(- \Delta)_p^s\, u & = b\, (u^+)^{p-1} - a\, (u^-)^{p-1} && \text{in } \Omega\\[10pt]
u & = 0 && \text{in } \R^N \setminus \Omega,
\end{aligned}\right.
\end{equation}
where $u^\pm = \max \set{\pm u,0}$ are the positive and negative parts of $u$, respectively, has a nontrivial weak solution. Let us recall the weak formulation of \eqref{1}. Let
\[
[u]_{s,p} = \left(\int_{\R^{2N}} \frac{|u(x) - u(y)|^p}{|x - y|^{N+sp}}\, dx dy\right)^{1/p}
\]
be the Gagliardo seminorm of the measurable function $u : \R^N \to \R$, and let
\[
W^{s,p}(\R^N) = \set{u \in L^p(\R^N) : [u]_{s,p} < \infty}
\]
be the fractional Sobolev space endowed with the norm
\[
\norm[s,p]{u} = \big(\pnorm[p]{u}^p + [u]_{s,p}^p\big)^{1/p},
\]
where $\pnorm[p]{\cdot}$ is the norm in $L^p(\R^N)$. We work in the closed linear subspace
\[
X_p^s(\Omega) = \set{u \in W^{s,p}(\R^N) : u = 0 \text{ a.e.\! in } \R^N \setminus \Omega}
\]
equivalently renormed by setting $\norm{\cdot} = [\cdot]_{s,p}$ (see Di Nezza et al.\! \cite[Theorem 7.1]{MR2944369}). A function $u \in X_p^s(\Omega)$ is a weak solution of problem \eqref{1} if
\begin{multline} \label{2}
\int_{\R^{2N}} \frac{|u(x) - u(y)|^{p-2}\, (u(x) - u(y))\, (v(x) - v(y))}{|x - y|^{N+sp}}\, dx dy\\[5pt]
= \int_\Omega \left[b\, (u^+)^{p-1} - a\, (u^-)^{p-1}\right] v\, dx \quad \forall v \in X_p^s(\Omega).
\end{multline}

This notion of spectrum for linear local elliptic partial differential operators was introduced by Dancer \cite{MR58:17506,MR82i:35063} and \Fucik \cite{MR56:5998}, who recognized its significance for the solvability of related semilinear boundary value problems. In particular, the Dancer-\Fucik spectrum of the Laplacian in $\Omega$ with the Dirichlet boundary condition is the set $\Sigma(\Omega)$ of all points $(a,b) \in \R^2$ such that the problem
\begin{equation} \label{3}
\left\{\begin{aligned}
- \Delta u & = b u^+ - a u^- && \text{in } \Omega\\[10pt]
u & = 0 && \text{on } \bdry{\Omega}
\end{aligned}\right.
\end{equation}
has a nontrivial solution. Denoting by $\lambda_k \nearrow + \infty$ the Dirichlet eigenvalues of $- \Delta$ in $\Omega$, $\Sigma(\Omega)$ clearly contains the sequence of points $(\lambda_k,\lambda_k)$. For $N = 1$, where $\Omega$ is an interval, \Fucik \cite{MR56:5998} showed that $\Sigma(\Omega)$ with the periodic boundary condition consists of a sequence of hyperbolic-like curves passing through the points $(\lambda_k,\lambda_k)$, with one or two curves going through each point. For $N \ge 2$, $\Sigma(\Omega)$ consists locally of curves emanating from the points $(\lambda_k,\lambda_k)$ (see Gallou{\"e}t and Kavian \cite{MR658734}, Ruf \cite{MR640779}, Lazer and McKenna \cite{MR871108}, Lazer \cite{MR965532}, C{\'a}c \cite{MR1011156}, Magalh{\~a}es \cite{MR1077275}, Cuesta and Gossez \cite{MR1181350}, de Figueiredo and Gossez \cite{MR1269657}, and Margulies and Margulies \cite{MR1484910}). Schechter \cite{MR1322614} showed that in the square $(\lambda_{k-1},\lambda_{k+1}) \times (\lambda_{k-1},\lambda_{k+1})$, $\Sigma(\Omega)$ contains two strictly decreasing curves, which may coincide, such that the points in the square that are either below the lower curve or above the upper curve are not in $\Sigma(\Omega)$, while the points between them may or may not belong to $\Sigma(\Omega)$ when they do not coincide.

The Dancer-\Fucik spectrum of the $p$-Laplacian $\Delta_p\, u = \divg \big(|\nabla u|^{p-2}\, \nabla u\big)$ is the set $\Sigma_p(\Omega)$ of all points $(a,b) \in \R^2$ such that the problem
\[
\left\{\begin{aligned}
- \Delta_p\, u & = b\, (u^+)^{p-1} - a\, (u^-)^{p-1} && \text{in } \Omega\\[10pt]
u & = 0 && \text{on } \bdry{\Omega}
\end{aligned}\right.
\]
has a nontrivial solution. For $N = 1$, the Dirichlet spectrum $\sigma(- \Delta_p)$ of $- \Delta_p$ in $\Omega$ consists of a sequence of simple eigenvalues $\lambda_k \nearrow + \infty$ and $\Sigma_p(\Omega)$ has the same general shape as $\Sigma(\Omega)$ (see Dr{\'a}bek \cite{MR94e:47084}). For $N \ge 2$, the first eigenvalue $\lambda_1$ of $- \Delta_p$ is positive, simple, and has an associated eigenfunction that is positive in $\Omega$ (see Anane \cite{MR89e:35124} and Lindqvist \cite{MR90h:35088,MR1139483}), so $\Sigma_p(\Omega)$ contains the two lines $\lambda_1 \times \R$ and $\R \times \lambda_1$. Moreover, $\lambda_1$ is isolated in the spectrum, so the second eigenvalue $\lambda_2 = \inf \sigma(- \Delta_p) \cap (\lambda_1,\infty)$ is well-defined (see Anane and Tsouli \cite{MR97k:35190}), and a first nontrivial curve in $\Sigma_p(\Omega)$ passing through $(\lambda_2,\lambda_2)$ and asymptotic to $\lambda_1 \times \R$ and $\R \times \lambda_1$ at infinity was constructed using the mountain pass theorem by Cuesta et al.\! \cite{MR1726923}. Although a complete description of $\sigma(- \Delta_p)$ is not yet available, an increasing and unbounded sequence of eigenvalues can be constructed via a standard minimax scheme based on the Krasnosel$'$ski\u\i\ genus, or via nonstandard schemes based on the cogenus as in Dr{\'a}bek and Robinson \cite{MR1726752} and the cohomological index as in Perera \cite{MR1998432}. Unbounded sequences of decreasing curves in $\Sigma_p(\Omega)$, analogous to the lower and upper curves of Schechter \cite{MR1322614} in the semilinear case, have been constructed using various minimax schemes by Cuesta \cite{Cuesta}, Micheletti and Pistoia \cite{MR2002d:35157}, and Perera \cite{MR2210289}.

Goyal and Sreenadh \cite{MR3238506} recently studied the Dancer-\Fucik spectrum for a class of linear nonlocal elliptic operators that includes the fractional Laplacian $(- \Delta)^s$. As in Cuesta et al.\! \cite{MR1726923}, they constructed a first nontrivial curve in the Dancer-\Fucik spectrum that passes through $(\lambda_2,\lambda_2)$ and is asymptotic to $\lambda_1 \times \R$ and $\R \times \lambda_1$ at infinity. Very recently, in \cite{brasco}, the authors proved, among other things, that the second variational eigenvalue $\lambda_2$ is larger than $\lambda_1$ and $(\lambda_1,\lambda_2)$ does not contain other eigevalues.

The purpose of this note is to point out that the general theories developed in Perera et al.\! \cite{MR2640827} and Perera and Schechter \cite{MR3012848} apply to the fractional $p$-Laplacian and Laplacian operators, respectively, and draw some conclusions about their Dancer-\Fucik spectra. We construct an unbounded sequence of decreasing curves in $\Sigma_p^s(\Omega)$ using a suitable minimax scheme. For $p = 2$, we present a very accurate local analysis. We construct the minimal and maximal curves of the spectrum locally near the points where it intersects the main diagonal of the plane. We give a sufficient condition for the region between them to be nonempty, and show that it is free of the spectrum in the case of a simple eigenvalue. Finally we compute the critical groups in various regions separated by these curves. We compute them precisely in certain regions, and prove a shifting theorem that gives a finite-dimensional reduction in certain other regions. This allows us to obtain nontrivial solutions of perturbed problems with nonlinearities crossing a curve of the spectrum via a comparison of the critical groups at zero and infinity.

\section{Dancer-\Fucik spectrum of the fractional $p$-Laplacian}

The general theory developed in Perera et al.\! \cite{MR2640827} applies to problem \eqref{1}. Indeed, the odd $(p - 1)$-homogeneous operator $A_p^s \in C(X_p^s(\Omega),X_p^s(\Omega)^\ast)$, where $X_p^s(\Omega)^\ast$ is the dual of $X_p^s(\Omega)$, defined by
\[
A_p^s(u)\, v = \int_{\R^{2N}} \frac{|u(x) - u(y)|^{p-2}\, (u(x) - u(y))\, (v(x) - v(y))}{|x - y|^{N+sp}}\, dx dy, \quad u, v \in X_p^s(\Omega)
\]
that is associated with the left-hand side of equation \eqref{2} satisfies
\begin{equation} \label{4}
A_p^s(u)\, u = \norm{u}^p, \quad |A_p^s(u)\, v| \le \norm{u}^{p-1} \norm{v} \quad \forall u, v \in X_p^s(\Omega)
\end{equation}
and is the Fr\'{e}chet derivative of the $C^1$-functional
\[
I_p^s(u) = \frac{1}{p} \int_{\R^{2N}} \frac{|u(x) - u(y)|^p}{|x - y|^{N+sp}}\, dx dy, \quad u \in X_p^s(\Omega).
\]
Moreover, since $X_p^s(\Omega)$ is uniformly convex, it follows from \eqref{4} that $A_p^s$ is of type (S), i.e., every sequence $\seq{u_j} \subset X_p^s(\Omega)$ such that
\[
u_j \wto u, \quad A_p^s(u_j)\, (u_j - u) \to 0
\]
has a subsequence that converges strongly to $u$ (see \cite[Proposition 1.3]{MR2640827}). Hence the operator $A_p^s$ satisfies the structural assumptions of \cite[Chapter 1]{MR2640827}.

When $a = b = \lambda$, \eqref{1} reduces to the nonlinear eigenvalue problem
\begin{equation} \label{5}
\left\{\begin{aligned}
(- \Delta)_p^s\, u & = \lambda\, |u|^{p-2}\, u && \text{in } \Omega\\[10pt]
u & = 0 && \text{in } \R^N \setminus \Omega.
\end{aligned}\right.
\end{equation}
Eigenvalues of this problem coincide with critical values of the functional
\[
\Psi(u) = \left(\int_\Omega |u|^p\, dx\right)^{-1}
\]
on the manifold
\[
\M = \set{u \in X_p^s(\Omega) : \norm{u} = 1}.
\]
The first eigenvalue
\[
\lambda_1 = \inf_{u \in \M}\, \Psi(u)
\]
is positive, simple, isolated, and has an associated eigenfunction that is positive in $\Omega$ (see Lindgren and Lindqvist \cite{MR3148135} and Franzina and Palatucci \cite {FrPa}), so $\Sigma_p^s(\Omega)$ contains the two lines $\lambda_1 \times \R$ and $\R \times \lambda_1$. Let $\F$ denote the class of symmetric subsets of $\M$, let $i(M)$ denote the $\Z_2$-cohomological index of $M \in \F$ (see Fadell and Rabinowitz \cite{MR57:17677}), and set
\[
\lambda_k := \inf_{\substack{M \in \F\\[1pt] i(M) \ge k}}\, \sup_{u \in M}\, \Psi(u), \quad k \ge 2.
\]
Then $\lambda_k \nearrow + \infty$ is a sequence of eigenvalues of problem \eqref{5} (see \cite[Theorem 4.6]{MR2640827}), so $\Sigma_p^s(\Omega)$ contains the sequence of points $(\lambda_k,\lambda_k)$.

Following \cite[Chapter 8]{MR2640827}, we now construct an unbounded sequence of decreasing curves in $\Sigma_p^s(\Omega)$. For $t > 0$, let
\[
\Psi_t(u) = \left(\int_\Omega \left[(u^+)^p + t\, (u^-)^p\right] dx\right)^{-1}, \quad u \in \M.
\]
Then the point $(c,c\, t) \in \Sigma_p^s(\Omega)$ if and only if $c$ is a critical value of $\Psi_t$ (see \cite[Lemma 8.3]{MR2640827}). For each $k \ge 2$ such that $\lambda_k > \lambda_{k-1}$, let
\[
C \Psi_t^{\lambda_{k-1}} = (\Psi_t^{\lambda_{k-1}} \times [0,1])/(\Psi_t^{\lambda_{k-1}} \times \set{1})
\]
be the cone on the sublevel set $\Psi_t^{\lambda_{k-1}} = \set{u \in \M : \Psi_t(u) \le \lambda_{k-1}}$, let $\Gamma_k$ denote the class of maps $\gamma \in C(C \Psi_t^{\lambda_{k-1}},\M)$ such that $\restr{\gamma}{\Psi_t^{\lambda_{k-1}}}$ is the identity, and set
\[
c_k^s(t) = \inf_{\gamma\in \Gamma_k}\, \sup_{u \in \gamma(C \Psi_t^{\lambda_{k-1}})}\, \Psi_t(u).
\]
We have the following theorem as a special case of \cite[Theorem 8.8]{MR2640827}.

\begin{theorem}
Let
\[
\C_k^s = \set{(c_k^s(t),c_k^s(t)\, t) : \lambda_{k-1}/\lambda_k < t < \lambda_k/\lambda_{k-1}}.
\]
Then $\C_k$ is a decreasing continuous curve in $\Sigma_p^s(\Omega)$, and $c_k^s(1) \ge \lambda_k$.
\end{theorem}

\section{Dancer-\Fucik spectrum of the fractional Laplacian}

The Dancer-\Fucik spectrum of the operator $(- \Delta)^s$ in $\Omega$ is the set $\Sigma^s(\Omega)$ of all points $(a,b) \in \R^2$ such that the problem
\begin{equation} \label{6}
\left\{\begin{aligned}
(- \Delta)^s\, u & = b u^+ - a u^- && \text{in } \Omega\\[10pt]
u & = 0 && \text{in } \R^N \setminus \Omega
\end{aligned}\right.
\end{equation}
has a nontrivial weak solution. The general theory developed in Perera and Schechter \cite{MR3012848} applies to problem \eqref{6}. Indeed, set $X^s(\Omega) = X_2^s(\Omega)$ and let $A^s$ be the inverse of the solution operator $S : L^2(\Omega) \to S(L^2(\Omega)) \subset X^s(\Omega),\, f \mapsto u$ of the problem
\[
\left\{\begin{aligned}
(- \Delta)^s\, u & = f(x) && \text{in } \Omega\\[10pt]
u & = 0 && \text{in } \R^N \setminus \Omega.
\end{aligned}\right.
\]
Then $A^s$ is a self-adjoint operator on $L^2(\Omega)$ and we have
\begin{align*}
& \ip{u}{v} = (A^{s/2} u,A^{s/2}v)_2=\int_{\R^{2N}} \frac{(u(x) - u(y))(v(x)-v(y)}{|x - y|^{N+2s}} dxdy, \quad \forall u, v \in X^s(\Omega), \\
& \norm{u} = \|A^{s/2} u\|_2=\Big(\int_{\R^{2N}} \frac{|u(x) - u(y)|^2}{|x - y|^{N+2s}}dxdy\Big)^{1/2}, \quad \forall u \in X^s(\Omega),
\end{align*}
where $\ip{\cdot}{\cdot}$ and $\ip[2]{\cdot}{\cdot}$ are the inner products in $X^s(\Omega)$ and $L^2(\Omega)$, respectively. Moreover, its spectrum $\sigma(A^s) \subset (0,\infty)$ and $(A^s)^{-1} : L^2(\Omega) \to L^2(\Omega)$ is a compact operator since the embedding $X^s(\Omega) \hookrightarrow L^2(\Omega)$ is compact. Therefore $\sigma(A^s)$ consists of isolated eigenvalues $\lambda_k,\, k \ge 1$ of finite multiplicities satisfying $0 < \lambda_1 < \lambda_2 < \cdots$. The first eigenvalue $\lambda_1$ is simple and has an associated eigenfunction $\varphi_1 > 0$, and if $w \in ((\R\, \varphi_1)^\perp \cap X^s(\Omega)) \setminus \set{0}$, then
\[
0 = \ip{w}{\varphi_1} = \ip[2]{A^s w}{\varphi_1} = \ip[2]{w}{A^s \varphi_1} = \lambda_1 \ip[2]{w}{\varphi_1},
\]
so $w^\pm \ne 0$. Hence the operator $A^s$ satisfies all the assumptions of \cite[Chapter 4]{MR3012848}.

Now we describe the minimal and maximal curves of $\Sigma^s(\Omega)$ in the square
\[
Q_k = (\lambda_{k-1},\lambda_{k+1})^2, \quad k \ge 2
\]
constructed in \cite{MR3012848}. Weak solutions of problem \eqref{6} coincide with critical points of the $C^1$-functional
\[
I(u,a,b) = \frac{1}{2} \int_{\R^{2N}} \frac{|u(x) - u(y)|^2}{|x - y|^{N+2s}}\, dx dy - \frac{1}{2} \int_\Omega \left[b\, (u^+)^2 + a\, (u^-)^2\right] dx, \quad u \in X^s(\Omega).
\]
Denote by $E_k$ the eigenspace of $\lambda_k$ and set
\[
N_k = \bigoplus_{j=1}^k E_j, \qquad M_k = N_k^\perp \cap X^s(\Omega).
\]
Then $X^s(\Omega) = N_k \oplus M_k$ is an orthogonal decomposition with respect to both $\ip{\cdot}{\cdot}$ and $\ip[2]{\cdot}{\cdot}$. When $(a,b) \in Q_k$, $I(v + y + w),\, v + y + w \in N_{k-1} \oplus E_k \oplus M_k$ is strictly concave in $v$ and strictly convex in $w$, i.e., if $v_1 \ne v_2 \in N_{k-1},\, w \in M_{k-1}$,
\[
I((1 - t)\, v_1 + t\, v_2 + w) > (1 - t)\, I(v_1 + w) + t\, I(v_2 + w) \quad \forall t \in (0,1),
\]
and if $v \in N_k,\, w_1 \ne w_2 \in M_k$,
\[
I(v + (1 - t)\, w_1 + t\, w_2) < (1 - t)\, I(v + w_1) + t\, I(v + w_2) \quad \forall t \in (0,1)
\]
(see \cite[Proposition 4.6.1]{MR3012848}).

\begin{proposition}[{\cite[Proposition 4.7.1, Corollary 4.7.3, Proposition 4.7.4]{MR3012848}}] \label{Proposition 1}
Let $(a,b) \in Q_k$.
\begin{enumroman}
\item There is a positive homogeneous map $\theta(\cdot,a,b) \in C(M_{k-1},N_{k-1})$ such that $v = \theta(w)$ is the unique solution of
    \[
    I(v + w) = \sup_{v' \in N_{k-1}} I(v' + w), \quad w \in M_{k-1}.
    \]
    Moreover, $I'(v + w) \perp N_{k-1}$ if and only if $v = \theta(w)$. Furthermore, $\theta$ is continuous on $M_{k-1} \times Q_k$ and $\theta(w,\lambda_k,\lambda_k) = 0$ for all $w \in M_{k-1}$.
\item There is a positive homogeneous map $\tau(\cdot,a,b) \in C(N_k,M_k)$ such that $w = \tau(v)$ is the unique solution of
    \[
    I(v + w) = \inf_{w' \in M_k} I(v + w'), \quad v \in N_k.
    \]
    Moreover, $I'(v + w) \perp M_k$ if and only if $w = \tau(v)$. Furthermore, $\tau$ is continuous on $N_k \times Q_k$ and $\tau(v,\lambda_k,\lambda_k) = 0$ for all $v \in N_k$.
\end{enumroman}
\end{proposition}

For $(a,b) \in Q_k$, let
\begin{gather*}
\sigma(w,a,b) = \theta(w,a,b) + w, \quad w \in M_{k-1}, \qquad S_k(a,b) = \sigma(M_{k-1},a,b),\\[10pt]
\zeta(v,a,b) = v + \tau(v,a,b), \quad v \in N_k, \qquad S^k(a,b) = \zeta(N_k,a,b).
\end{gather*}
Then $S_k$ and $S^k$ are topological manifolds modeled on $M_{k-1}$ and $N_k$, respectively. Thus, $S_k$ is infinite dimensional, while $S^k$ is $d_k$-dimensional, where $d_k = \dim N_k$. For $B \subset X^s(\Omega)$, set $\widetilde{B} = \bgset{u \in B : \norm{u} = 1}$. We say that $B$ is a radial set if $B = \bgset{tu : u \in \widetilde{B},\, t \ge 0}$. Since $\theta$ and $\tau$ are positive homogeneous, so are $\sigma$ and $\zeta$, and hence $S_k$ and $S^k$ are radial manifolds.

Let
\[
K(a,b) = \bgset{u \in X^s(\Omega) : I'(u,a,b) = 0}
\]
be the set of critical points of $I(\cdot,a,b)$. Since $I'$ is positive homogeneous, $K$ is a radial set. Since $I(u) = \ip{I'(u)}{u}\!/2$,
\begin{equation} \label{7}
I(u) = 0 \quad \forall u \in K.
\end{equation}
Since $X^s(\Omega) = N_{k-1} \oplus E_k \oplus M_k$, Proposition \ref{Proposition 1} implies
\begin{equation} \label{8}
K = \bgset{u \in S_k \cap S^k : I'(u) \perp E_k}.
\end{equation}
Together with \eqref{7}, it also implies
\begin{equation} \label{9}
K \subset \bgset{u \in S_k \cap S^k : I(u) = 0}.
\end{equation}

Set
\begin{gather*}
n_{k-1}(a,b) = \inf_{w \in \widetilde{M}_{k-1}} \sup_{v \in N_{k-1}} I(v + w,a,b),\\[10pt]
m_k(a,b) = \sup_{v \in \widetilde{N}_k} \inf_{w \in M_k} I(v + w,a,b).
\end{gather*}
Since $I(u,a,b)$ is nonincreasing in $a$ for fixed $u$ and $b$, and in $b$ for fixed $u$ and $a$, $n_{k-1}(a,b)$ and $m_k(a,b)$ are nonincreasing in $a$ for fixed $b$, and in $b$ for fixed $a$. By Proposition \ref{Proposition 1},
\begin{gather*}
n_{k-1}(a,b) = \inf_{w \in \widetilde{M}_{k-1}} I(\sigma(w,a,b),a,b),\\[10pt]
m_k(a,b) = \sup_{v \in \widetilde{N}_k} I(\zeta(v,a,b),a,b).
\end{gather*}

\begin{proposition}[{\cite[Proposition 4.7.5, Lemma 4.7.6, Proposition 4.7.7]{MR3012848}}]
Let $(a,b), (a',b') \in Q_k$.
\begin{enumroman}
\item Assume that $n_{k-1}(a,b) = 0$. Then
    \begin{gather*}
    I(u,a,b) \ge 0 \quad \forall u \in S_k(a,b),\\[10pt]
    K(a,b) = \bgset{u \in S_k(a,b) : I(u,a,b) = 0},
    \end{gather*}
    and $(a,b) \in \Sigma^s(\Omega)$.
    \begin{enumalph}
    \item If $a' \le a$, $b' \le b$, and $(a',b') \ne (a,b)$, then $n_{k-1}(a',b') > 0$,
        \[
        I(u,a',b') > 0 \quad \forall u \in S_k(a',b') \setminus \set{0},
        \]
        and $(a',b') \notin \Sigma^s(\Omega)$.
    \item If $a' \ge a$, $b' \ge b$, and $(a',b') \ne (a,b)$, then $n_{k-1}(a',b') < 0$ and there is a $u \in S_k(a',b') \setminus \set{0}$ such that
        \[
        I(u,a',b') < 0.
        \]
    \end{enumalph}
    Furthermore, $n_{k-1}$ is continuous on $Q_k$ and $n_{k-1}(\lambda_k,\lambda_k) = 0$.
\item Assume that $m_k(a,b) = 0$. Then
    \begin{gather*}
    I(u,a,b) \le 0 \quad \forall u \in S^k(a,b),\\[10pt]
    K(a,b) = \bgset{u \in S^k(a,b) : I(u,a,b) = 0},
    \end{gather*}
    and $(a,b) \in \Sigma^s(\Omega)$.
    \begin{enumalph}
    \item If $a' \ge a$, $b' \ge b$, and $(a',b') \ne (a,b)$, then $m_k(a',b') < 0$,
        \[
        I(u,a',b') < 0 \quad \forall u \in S^k(a',b') \setminus \set{0},
        \]
        and $(a',b') \notin \Sigma^s(\Omega)$.
    \item If $a' \le a$, $b' \le b$, and $(a',b') \ne (a,b)$, then $m_k(a',b') > 0$ and there is a $u \in S^k(a',b') \setminus \set{0}$ such that
        \[
        I(u,a',b') > 0.
        \]
    \end{enumalph}
    Furthermore, $m_k$ is continuous on $Q_k$ and $m_k(\lambda_k,\lambda_k) = 0$.
\end{enumroman}
\end{proposition}

For $a \in (\lambda_{k-1},\lambda_{k+1})$, set
\begin{gather*}
\nu_{k-1}(a) = \sup \bgset{b \in (\lambda_{k-1},\lambda_{k+1}) : n_{k-1}(a,b) \ge 0},\\[10pt]
\mu_k(a) = \inf \bgset{b \in (\lambda_{k-1},\lambda_{k+1}) : m_k(a,b) \le 0}.
\end{gather*}
Then
\begin{gather*}
b = \nu_{k-1}(a) \iff n_{k-1}(a,b) = 0,\\[10pt]
b = \mu_k(a) \iff m_k(a,b) = 0
\end{gather*}
(see \cite[Lemma 4.7.8]{MR3012848}).

\begin{theorem}[{\cite[Theorem 4.7.9]{MR3012848}}]
Let $(a,b) \in Q_k$.
\begin{enumroman}
\item The function $\nu_{k-1}$ is continuous, strictly decreasing, and satisfies
    \begin{enumalph}
    \item $\nu_{k-1}(\lambda_k) = \lambda_k$,
    \item $b = \nu_{k-1}(a) \implies (a,b) \in \Sigma^s(\Omega)$,
    \item $b < \nu_{k-1}(a) \implies (a,b) \notin \Sigma^s(\Omega)$.
    \end{enumalph}
\item The function $\mu_k$ is continuous, strictly decreasing, and satisfies
    \begin{enumalph}
    \item $\mu_k(\lambda_k) = \lambda_k$,
    \item $b = \mu_k(a) \implies (a,b) \in \Sigma^s(\Omega)$,
    \item $b > \mu_k(a) \implies (a,b) \notin \Sigma^s(\Omega)$.
    \end{enumalph}
\item $\nu_{k-1}(a) \le \mu_k(a)$.
\end{enumroman}
\end{theorem}

Thus,
\[
C_k : b = \nu_{k-1}(a), \qquad C^k : b = \mu_k(a)
\]
are strictly decreasing curves in $Q_k$ that belong to $\Sigma^s(\Omega)$. They both pass through the point $(\lambda_k,\lambda_k)$ and may coincide. The region
\[
\Il{k} = \bgset{(a,b) \in Q_k : b < \nu_{k-1}(a)}
\]
below the lower curve $C_k$ and the region
\[
\Iu{k} = \bgset{(a,b) \in Q_k : b > \mu_k(a)}
\]
above the upper curve $C^k$ are free of $\Sigma^s(\Omega)$. They are the minimal and maximal curves of $\Sigma^s(\Omega)$ in $Q_k$ in this sense. Points in the region
\[
\II{k} = \bgset{(a,b) \in Q_k : \nu_{k-1}(a) < b < \mu_k(a)}
\]
between $C_k$ and $C^k$, when it is nonempty, may or may not belong to $\Sigma^s(\Omega)$.

For $(a,b) \in Q_k$, let
\[
\calN_k(a,b) = S_k(a,b) \cap S^k(a,b).
\]
Since $S_k$ and $S^k$ are radial sets, so is $\calN_k$. The next two propositions show that $\calN_k$ is a topological manifold modeled on $E_k$ and hence
\[
\dim \calN_k = d_k - d_{k-1}.
\]
We will call it the null manifold of $I$.

\begin{proposition}[{\cite[Proposition 4.8.1, Lemma 4.8.3, Proposition 4.8.4]{MR3012848}}]
Let $(a,b) \in Q_k$.
\begin{enumroman}
\item There is a positive homogeneous map $\eta(\cdot,a,b) \in C(E_k,N_{k-1})$ such that $v = \eta(y)$ is the unique solution of
    \[
    I(\zeta(v + y)) = \sup_{v' \in N_{k-1}} I(\zeta(v' + y)), \quad y \in E_k.
    \]
    Moreover, $I'(\zeta(v + y)) \perp N_{k-1}$ if and only if $v = \eta(y)$. Furthermore, $\eta$ is continuous on $E_k \times Q_k$ and $\eta(y,\lambda_k,\lambda_k) = 0$ for all $y \in E_k$.
\item There is a positive homogeneous map $\xi(\cdot,a,b) \in C(E_k,M_k)$ such that $w = \xi(y)$ is the unique solution of
    \[
    I(\sigma(y + w)) = \inf_{w' \in M_k} I(\sigma(y + w')), \quad y \in E_k.
    \]
    Moreover, $I'(\sigma(y + w)) \perp M_k$ if and only if $w = \xi(y)$. Furthermore, $\xi$ is continuous on $E_k \times Q_k$ and $\xi(y,\lambda_k,\lambda_k) = 0$ for all $y \in E_k$.
\item For all $y \in E_k$, $\zeta(\eta(y) + y) = \sigma(y + \xi(y))$, i.e., $\eta(y) = \theta(y + \xi(y))$ and $\xi(y) = \tau(\eta(y) + y)$.
\end{enumroman}
\end{proposition}

Let
\[
\varphi(y) = \zeta(\eta(y) + y) = \sigma(y + \xi(y)), \quad y \in E_k.
\]

\begin{proposition}[{\cite[Proposition 4.8.5]{MR3012848}}]
Let $(a,b) \in Q_k$.
\begin{enumroman}
\item $\varphi(\cdot,a,b) \in C(E_k,X^s(\Omega))$ is a positive homogeneous map such that
\[
I(\varphi(y)) = \inf_{w \in M_k} \sup_{v \in N_{k-1}} I(v + y + w) = \sup_{v \in N_{k-1}} \inf_{w \in M_k} I(v + y + w), \quad y \in E_k
\]
and $I'(\varphi(y)) \in E_k$ for all $y \in E_k$.
\item If $(a',b') \in Q_k$ with $a' \ge a$ and $b' \ge b$, then
\[
I(\varphi(y,a',b'),a',b') \le I(\varphi(y,a,b),a,b) \quad \forall y \in E_k.
\]
\item $\varphi$ is continuous on $E_k \times Q_k$.
\item $\varphi(y,\lambda_k,\lambda_k) = y \quad \forall y \in E_k$.
\item $\calN_k(a,b) = \bgset{\varphi(y,a,b) : y \in E_k}$.
\item $\calN_k(\lambda_k,\lambda_k) = E_k$.
\end{enumroman}
\end{proposition}

By \eqref{8} and \eqref{9},
\begin{equation} \label{10}
K = \bgset{u \in \calN_k : I'(u) \perp E_k} \subset \bgset{u \in \calN_k : I(u) = 0}.
\end{equation}
The following theorem shows that the curves $C_k$ and $C^k$ are closely related to $\widetilde{I} = \restr{I}{\calN_k}$.

\begin{theorem}[{\cite[Theorem 4.8.6]{MR3012848}}]
Let $(a,b) \in Q_k$.
\begin{enumroman}
\item If $b < \nu_{k-1}(a)$, then
    \[
    \widetilde{I}(u,a,b) > 0 \quad \forall u \in \calN_k(a,b) \setminus \set{0}.
    \]
\item If $b = \nu_{k-1}(a)$, then
    \begin{gather*}
    \widetilde{I}(u,a,b) \ge 0 \quad \forall u \in \calN_k(a,b),\\[10pt]
    K(a,b) = \bgset{u \in \calN_k(a,b) : \widetilde{I}(u,a,b) = 0}.
    \end{gather*}
\item If $\nu_{k-1}(a) < b < \mu_k(a)$, then there are $u_i \in \calN_k(a,b) \setminus \set{0},\, i = 1, 2$ such that
    \[
    \widetilde{I}(u_1,a,b) < 0 < \widetilde{I}(u_2,a,b).
    \]
\item If $b = \mu_k(a)$, then
    \begin{gather*}
    \widetilde{I}(u,a,b) \le 0 \quad \forall u \in \calN_k(a,b),\\[10pt]
    K(a,b) = \bgset{u \in \calN_k(a,b) : \widetilde{I}(u,a,b) = 0}.
    \end{gather*}
\item If $b > \mu_k(a)$, then
    \[
    \widetilde{I}(u,a,b) < 0 \quad \forall u \in \calN_k(a,b) \setminus \set{0}.
    \]
\end{enumroman}
\end{theorem}

By \eqref{10}, solutions of \eqref{6} are in $\calN_k$. The set $K(a,b)$ of solutions is all of $\calN_k(a,b)$ exactly when $(a,b) \in Q_k$ is on both $C_k$ and $C^k$ (see \cite[Theorem 4.8.7]{MR3012848}). When $\lambda_k$ is a simple eigenvalue, $\calN_k$ is $1$-dimensional and hence this implies that $(a,b)$ is on exactly one of those curves if and only if
\[
K(a,b) = \bgset{t\, \varphi(y_0,a,b) : t \ge 0}
\]
for some $y_0 \in E_k \setminus \set{0}$ (see \cite[Corollary 4.8.8]{MR3012848}).

The following theorem gives a sufficient condition for the region $\II{k}$ to be nonempty.

\begin{theorem}[{\cite[Theorem 4.9.1]{MR3012848}}]
If there is a $y \in E_k$ such that $\spnorm[2]{y^+} \ne \spnorm[2]{y^-}$, then there is a neighborhood $N \subset Q_k$ of $(\lambda_k,\lambda_k)$ such that every point $(a,b) \in N \setminus \set{(\lambda_k,\lambda_k)}$ with $a + b = 2 \lambda_k$ is in {\em $\II{k}$}.
\end{theorem}

For the local problem \eqref{3}, this result is due to Li et al. \cite{MR2389991}. When $\lambda_k$ is a simple eigenvalue, the region $\II{k}$ is free of $\Sigma^s(\Omega)$ (see \cite[Theorem 4.10.1]{MR3012848}). For problem \eqref{3}, this is due to Gallou{\"e}t and Kavian \cite{MR658734}.

When $(a,b) \notin \Sigma^s(\Omega)$, $0$ is the only critical point of $I$ and its critical groups are given by
\[
C_q(I,0) = H_q(I^0,I^0 \setminus \set{0}), \quad q \ge 0,
\]
where $I^0 = \bgset{u \in X^s(\Omega) : I(u) \le 0}$ and $H$ denotes singular homology. We take the coefficient group to be the field $\Z_2$. The following theorem gives our main results on the critical groups.

\begin{theorem}[{\cite[Theorem 4.11.2]{MR3012848}}]
Let $(a,b) \in Q_k \setminus \Sigma^s(\Omega)$.
\begin{enumroman}
\item If {\em $(a,b) \in \Il{k}$}, then
\[
C_q(I,0) \isom \delta_{qd_{k-1}}\, \Z_2.
\]
\item If {\em $(a,b) \in \Iu{k}$}, then
\[
C_q(I,0) \isom \delta_{qd_k}\, \Z_2.
\]
\item If {\em $(a,b) \in \II{k}$}, then
\[
C_q(I,0) = 0, \quad q \le d_{k-1} \text{ or } q \ge d_k
\]
and
\[
C_q(I,0) \isom \widetilde{H}_{q - d_{k-1} - 1}(\set{u \in \calN_k : I(u) < 0}), \quad d_{k-1} < q < d_k,
\]
where $\widetilde{H}$ denotes the reduced homology groups. In particular, $C_q(I,0) = 0$ for all $q$ when $\lambda_k$ is simple.
\end{enumroman}
\end{theorem}

For the local problem \eqref{3}, this result is due to Dancer \cite{MR1725568,MR1921046} and Perera and Schechter \cite{MR1701261,MR1807945,MR1823913}. It can be used, for example, to obtain nontrivial solutions of perturbed problems with nonlinearities that cross a curve of the Dancer-\Fucik spectrum, via a comparison of the critical groups at zero and infinity. Consider the problem
\begin{equation} \label{12}
\left\{\begin{aligned}
(- \Delta)^s\, u & = f(x,u) && \text{in } \Omega\\[10pt]
u & = 0 && \text{in } \R^N \setminus \Omega,
\end{aligned}\right.
\end{equation}
where $f$ is a Carath\'{e}odory function on $\Omega \times \R$.

\begin{theorem}[{see \cite[Theorem 5.6.1]{MR3012848}}]
If
\[
f(x,t) = \begin{cases}
b_0 t^+ - a_0 t^- + \o(t) & \text{as } t \to 0\\[10pt]
b t^+ - a t^- + \o(t) & \text{as } |t| \to \infty,
\end{cases}
\]
uniformly a.e.\! in $\Omega$, for some $(a_0,b_0)$ and $(a,b)$ in $Q_k \setminus \Sigma^s(\Omega)$ that are on opposite sides of $C_k$ or $C^k$, then problem \eqref{12} has a nontrivial weak solution.
\end{theorem}

For problem \eqref{3}, this was proved in Perera and Schechter \cite{MR1807945}. It generalizes a well\nobreakdash-known result of Amann and Zehnder \cite{MR82b:47077} on the existence of nontrivial solutions for problems crossing an eigenvalue.

\def\cdprime{$''$}


\begin{thebibliography}{10}

\bibitem{MR82b:47077}
H.~Amann and E.~Zehnder.
\newblock Nontrivial solutions for a class of nonresonance problems and
  applications to nonlinear differential equations.
\newblock {\em Ann. Scuola Norm. Sup. Pisa Cl. Sci. (4)}, 7(4):539--603, 1980.

\bibitem{MR97k:35190}
A.~Anane and N.~Tsouli.
\newblock On the second eigenvalue of the {$p$}-{L}aplacian.
\newblock In {\em Nonlinear partial differential equations (F\`es, 1994)},
  volume 343 of {\em Pitman Res. Notes Math. Ser.}, pages 1--9. Longman,
  Harlow, 1996.

\bibitem{brasco}
L.\ Brasco and E.\ Parini.
\newblock The second eigenvalue of the fractional $p$-Laplacian.
\newblock Preprint.

\bibitem{MR89e:35124}
A.~Anane.
\newblock Simplicit\'e et isolation de la premi\`ere valeur propre du
  {$p$}-laplacien avec poids.
\newblock {\em C. R. Acad. Sci. Paris S\'er. I Math.}, 305(16):725--728, 1987.

\bibitem{MR1011156}
N.P.\ C{\'a}c.
\newblock On nontrivial solutions of a {D}irichlet problem whose jumping
  nonlinearity crosses a multiple eigenvalue.
\newblock {\em J. Differential Equations}, 80(2):379--404, 1989.

\bibitem{Ca}
L.\ Caffarelli.
\newblock Non-local diffusions, drifts and games.
\newblock In {\em Nonlinear Partial Differential Equations}, volume~7 of {\em
  Abel Symposia}, pages 37--52, 2012.

\bibitem{Cuesta}
M.~Cuesta.
\newblock On the {F}u\v c\'\i k spectrum of the {L}aplacian and the
  {$p$}-{L}aplacian.
\newblock In {\em Proceedings of Seminar in Differential Equations}, pages
  67--96, Kvilda, Czech Republic, May 29 -- June 2, 2000. Centre of Applied
  Mathematics, Faculty of Applied Sciences, University of West Bohemia in
  Pilsen.

\bibitem{MR1726923}
M.~Cuesta, D.~de~Figueiredo, and J.-P. Gossez.
\newblock The beginning of the {F}u\v cik spectrum for the {$p$}-{L}aplacian.
\newblock {\em J. Differential Equations}, 159(1):212--238, 1999.

\bibitem{MR1181350}
M.\ Cuesta and J.-P.\ Gossez.
\newblock A variational approach to nonresonance with respect to the {F}u\v cik
  spectrum.
\newblock {\em Nonlinear Anal.}, 19(5):487--500, 1992.

\bibitem{MR58:17506}
E.N. Dancer.
\newblock On the {D}irichlet problem for weakly non-linear elliptic partial
  differential equations.
\newblock {\em Proc. Roy. Soc. Edinburgh Sect. A}, 76(4):283--300, 1976/77.

\bibitem{MR82i:35063}
E.N. Dancer.
\newblock Corrigendum: ``{O}n the {D}irichlet problem for weakly nonlinear
  elliptic partial differential equations'' [{P}roc. {R}oy. {S}oc. {E}dinburgh
  {S}ect. {A} {\bf 76} (1976/77), no. 4, 283--300;\ {MR} {\bf 58} \#17506].
\newblock {\em Proc. Roy. Soc. Edinburgh Sect. A}, 89(1-2):15, 1981.

\bibitem{MR1725568}
E.N. Dancer.
\newblock Remarks on jumping nonlinearities.
\newblock In {\em Topics in nonlinear analysis}, volume~35 of {\em Progr.
  Nonlinear Differential Equations Appl.}, pages 101--116. Birkh\"auser, Basel,
  1999.

\bibitem{MR1921046}
E.N. Dancer.
\newblock Some results for jumping nonlinearities.
\newblock {\em Topol. Methods Nonlinear Anal.}, 19(2):221--235, 2002.

\bibitem{MR1269657}
D.G. de~Figueiredo and J.-P. Gossez.
\newblock On the first curve of the {F}u\v cik spectrum of an elliptic
  operator.
\newblock {\em Differential Integral Equations}, 7(5-6):1285--1302, 1994.

\bibitem{DiKuPa}
A.~Di~Castro, T.~Kuusi, and G.~Palatucci.
\newblock Local behavior of fractional $p$-minimizers.
\newblock preprint.

\bibitem{MR2944369}
E.\ Di~Nezza, G.\ Palatucci, and E.\ Valdinoci.
\newblock Hitchhiker's guide to the fractional {S}obolev spaces.
\newblock {\em Bull. Sci. Math.}, 136(5):521--573, 2012.

\bibitem{MR94e:47084}
P.~Dr{\'a}bek.
\newblock {\em Solvability and bifurcations of nonlinear equations}, volume 264
  of {\em Pitman Research Notes in Mathematics Series}.
\newblock Longman Scientific \& Technical, Harlow, 1992.

\bibitem{MR1726752}
P.\ Dr{\'a}bek and S.B. Robinson.
\newblock Resonance problems for the {$p$}-{L}aplacian.
\newblock {\em J. Funct. Anal.}, 169(1):189--200, 1999.

\bibitem{MR57:17677}
E.R. Fadell and P.H. Rabinowitz.
\newblock Generalized cohomological index theories for {L}ie group actions with
  an application to bifurcation questions for {H}amiltonian systems.
\newblock {\em Invent. Math.}, 45(2):139--174, 1978.

\bibitem{FrPa}
G.~Franzina and G.~Palatucci.
\newblock Fractional $p$-eigenvalues.
\newblock{\em Riv. Mat. Univ. Parma}, to appear.

\bibitem{MR56:5998}
S.\ Fu{\v{c}}{\'{\i}}k.
\newblock Boundary value problems with jumping nonlinearities.
\newblock {\em \v Casopis P\v est. Mat.}, 101(1):69--87, 1976.

\bibitem{MR658734}
T.\ Gallou{\"e}t and O.\ Kavian.
\newblock R\'esultats d'existence et de non-existence pour certains probl\`emes
  demi-lin\'eaires \`a l'infini.
\newblock {\em Ann. Fac. Sci. Toulouse Math. (5)}, 3(3-4):201--246 (1982),
  1981.

\bibitem{MR3238506}
Sarika Goyal and K.~Sreenadh.
\newblock On the {F}u\v cik spectrum of non-local elliptic operators.
\newblock {\em NoDEA Nonlinear Differential Equations Appl.}, 21(4):567--588,
  2014.

\bibitem{IaLiPeSq}
A.~Iannizzotto, S.~Liu, K.~Perera, and M.~Squassina.
\newblock Existence results for fractional $p$-{L}aplacian problems via {M}orse
  theory.
\newblock preprint.

\bibitem{MR3245079}
A.\ Iannizzotto and M.\ Squassina.
\newblock Weyl-type laws for fractional {$p$}-eigenvalue problems.
\newblock {\em Asymptot. Anal.}, 88(4):233--245, 2014.

\bibitem{MR871108}
A.C. Lazer and P.J. McKenna.
\newblock Critical point theory and boundary value problems with nonlinearities
  crossing multiple eigenvalues. {II}.
\newblock {\em Comm. Partial Differential Equations}, 11(15):1653--1676, 1986.

\bibitem{MR965532}
A.\ Lazer.
\newblock Introduction to multiplicity theory for boundary value problems with
  asymmetric nonlinearities.
\newblock In {\em Partial differential equations (Rio de Janeiro, 1986)},
  volume 1324 of {\em Lecture Notes in Math.}, pages 137--165. Springer,
  Berlin, 1988.

\bibitem{MR2389991}
C.\ Li, S.\ Li, and Z.\ Liu.
\newblock Existence of type ({II}) regions and convexity and concavity of
  potential functionals corresponding to jumping nonlinear problems.
\newblock {\em Calc. Var. Partial Differential Equations}, 32(2):237--251,
  2008.

\bibitem{MR3148135}
E.\ Lindgren and P.\ Lindqvist.
\newblock Fractional eigenvalues.
\newblock {\em Calc. Var. Partial Differential Equations}, 49(1-2):795--826,
  2014.

\bibitem{MR90h:35088}
P.\ Lindqvist.
\newblock On the equation {${\rm div}\,(\vert \nabla u\vert \sp {p-2}\nabla
  u)+\lambda\vert u\vert \sp {p-2}u=0$}.
\newblock {\em Proc. Amer. Math. Soc.}, 109(1):157--164, 1990.

\bibitem{MR1139483}
P.\ Lindqvist.
\newblock Addendum: ``{O}n the equation {${\rm div}(\vert \nabla u\vert \sp
  {p-2}\nabla u)+\lambda\vert u\vert \sp {p-2}u=0$}'' [{P}roc.\ {A}mer.\
  {M}ath.\ {S}oc.\ 109 (1990), no.\ 1, 157--164; {MR} 90h:35088].
\newblock {\em Proc. Amer. Math. Soc.}, 116(2):583--584, 1992.

\bibitem{MR1077275}
C.A. Magalh{\~a}es.
\newblock Semilinear elliptic problem with crossing of multiple eigenvalues.
\newblock {\em Comm. Partial Differential Equations}, 15(9):1265--1292, 1990.

\bibitem{MR1484910}
C.A.\ Margulies and W.\ Margulies.
\newblock An example of the {F}u\v cik spectrum.
\newblock {\em Nonlinear Anal.}, 29(12):1373--1378, 1997.

\bibitem{MR2002d:35157}
A.M.\ Micheletti and A.\ Pistoia.
\newblock On the {F}u\v c\'\i k spectrum for the {$p$}-{L}aplacian.
\newblock {\em Differential Integral Equations}, 14(7):867--882, 2001.

\bibitem{MR1998432}
K.\ Perera.
\newblock Nontrivial critical groups in {$p$}-{L}aplacian problems via the
  {Y}ang index.
\newblock {\em Topol. Methods Nonlinear Anal.}, 21(2):301--309, 2003.

\bibitem{MR2210289}
K.\ Perera.
\newblock On the {F}u\v c\'\i k spectrum of the {$p$}-{L}aplacian.
\newblock {\em NoDEA Nonlinear Differential Equations Appl.}, 11(2):259--270,
  2004.

\bibitem{MR2640827}
K.\ Perera, R.P.\ Agarwal, and D.\ O'Regan.
\newblock {\em Morse theoretic aspects of {$p$}-{L}aplacian type operators},
  volume 161 of {\em Mathematical Surveys and Monographs}.
\newblock American Mathematical Society, Providence, RI, 2010.

\bibitem{MR1701261}
K.\ Perera and M.\ Schechter.
\newblock Type {II} regions between curves of the {F}u\v cik spectrum and
  critical groups.
\newblock {\em Topol. Methods Nonlinear Anal.}, 12(2):227--243, 1998.

\bibitem{MR1807945}
K.\ Perera and M.\ Schechter.
\newblock A generalization of the {A}mann-{Z}ehnder theorem to nonresonance
  problems with jumping nonlinearities.
\newblock {\em NoDEA Nonlinear Differential Equations Appl.}, 7(4):361--367,
  2000.

\bibitem{MR1823913}
K.\ Perera and M.\ Schechter.
\newblock The {F}u\v c\'\i k spectrum and critical groups.
\newblock {\em Proc. Amer. Math. Soc.}, 129(8):2301--2308 (electronic), 2001.

\bibitem{MR3012848}
K.\ Perera and M.\ Schechter.
\newblock {\em Topics in critical point theory}, volume 198 of {\em Cambridge
  Tracts in Mathematics}.
\newblock Cambridge University Press, Cambridge, 2013.

\bibitem{MR640779}
B.\ Ruf.
\newblock On nonlinear elliptic problems with jumping nonlinearities.
\newblock {\em Ann. Mat. Pura Appl. (4)}, 128:133--151, 1981.

\bibitem{MR1322614}
M.\ Schechter.
\newblock The {F}u\v c\'\i k spectrum.
\newblock {\em Indiana Univ. Math. J.}, 43(4):1139--1157, 1994.

\end{thebibliography}
\end{document}